\documentclass[reqno]{amsart}
\usepackage{amsmath,amsfonts,amssymb,amsthm,mathrsfs,bbm,stmaryrd,array,url}

\makeatletter\def\@biblabel#1{#1}\makeatother
\newtheorem{theorem}{Theorem}[section]
\newtheorem{corollary}{Corollary}[section]
\newtheorem{lemma}{Lemma}[section]

\newtheorem*{conj1}{The Prime Twins Conjecture}

\newtheorem{definition}{Definition}[section]

\def\mydash{\CJKglue\raise0.2ex\hbox{---\kern-0.01em---}\CJKglue}

\begin{document}
\title[Proof of infinitely many pairs of primes $p$ and $p+2$ by induction to absurdity]{Proof of infinitely many pairs of primes $p$ and $p+2$ by induction to absurdity$^*$\footnotetext{*~An appendix is added at the end of the paper.}}
\author[Guangchang Dong]{Guangchang Dong}
\address{\newline\indent Guangchang Dong \newline\indent Department of mathematics\newline\indent Zhejiang University, Hangzhou, Zhejiang, China}
\email[Guangchang Dong]{dogc@zju.edu.cn}

\begin{abstract}
We take the pre-sieved set to be all natural numbers $N=\{1,2,3,\dots\}$ with a sieve system: single sieve, double sieve, $\cdots$

With single sieve, i.e., remove out the multiple of a prime, we derive all the primes. With double sieve, i.e., remove out the multiple and the multiple of a prime and $-2$ simultaneously, we get all the prime twins and prove that infinitely many prime twins exist under suitable conditions. Finally, with special 4 sieve, we prove that infinitely many prime twins exist without any restriction.
\end{abstract}

\maketitle
\tableofcontents

\section{Introduction}

\begin{conj1}[Euclid about 300 B.C.]
  There exist infinite many prime pairs $p$ and $p+2$, such as $3,5 ; 5,7 ; 11,13 ; 17,19 ; \cdots$
\end{conj1}
Chen[1] made a remarkable effort on this problem, he proved that, there exist infinite many $p$ and $p+2$ such that $p$ is a prime, $p+2$ is the almost prime (at most two prime factors).

Recently, Zhang[2] improved a partial result of Chen to
\begin{equation*}
  \lim \inf_{n\rightarrow\infty}(p_{n+1}-p_{n})<7\times10^{7},
\end{equation*}
where $p_{n}$ is the $n$'th prime.

However, the real conjecture has not solved yet.

Chen[1] used analytic theory of numbers. This theory based on sieve system collected in [3], i.e., Eratosthenes sieve, Brun's sieve, Rosser's sieve, Selberg's upper sieve, .... All those sieve methods have a common point, that is, their pre-sieve set are finite. It is obvious that sieve method with finite pre-sieve set cannot derive foundational property of primes: there exist infinitely many primes.

Since the infinitely many primes is a corollary of prime twins conjecture, so that the sieve system of [3] can not prove the prime twins conjecture. In order to prove the prime twins conjecture, a new sieve system is needed.

Take the pre-sieve set to be all natural numbers $N=\{1,2,3,4,5,...\}$. The sieve sets are determined progressively.

The ascending period of $N$ is 1, the smallest number except 1 is the prime $2=p_{0}$.

Denote $S_{p}$ be the sieve operator that remove all the multiple of p.

\textbf{Sieve step 0}. Apply $S_{p_{0}}$ to $N$, we have
\begin{equation*}
S_{p_{0}}N\equiv N_{0}=\{1,3,5,7,...\} .
\end{equation*}

We call $N_{0}$ the post-sieve set. The period of $N_{0}$ is 2, the number of elements in a period is 1, the smallest number except 1 is the prime $3=p_{1}$.

\textbf{Sieve step 1}.
\begin{equation*}
S_{p_{1}}N_{0}\equiv N_{1}=\{1,5,7,11,13,17,19,23,25,29,31,...\}
\end{equation*}

The ascending period of post-sieve set $N_{1}$ is $6=p_{0}p_{1}$, the number of elements in a period is $2=(p_{0}-1)(p_{1}-1)$, the smallest number except 1 is the prime $5=p_{2}$.

\textbf{Sieve step 2}.
\begin{equation*}
S_{p_{2}}N_{1}\equiv N_{2}=\{1,7,11,13,17,19,23,29,31,...\}
\end{equation*}

The ascending period of $N_{2}$ is $30=p_{0}p_{1}p_{2}$, the number of elements in a period is $8=(p_{0}-1)(p_{1}-1)(p_{2}-1)$, the smallest number except 1 is the prime $7=p_{3}$.

\textbf{Sieve step 3}.Remove the multiple of 7 we obtain $N_{3},...$.

From the above procedure we obtain the following general rule.

Perform the m's sieve step we obtain $N_{m}$. The ascending period of $N_{m}$ is $\prod\limits_{i=0}\limits^{m} p_{i}$, the number of elements in a period is $\prod\limits_{i=0}\limits^{m} (p_{i}-1)$. Denote the smallest number except 1 by $p_{m+1}$, then $p_{m+1}$ must be a prime. Otherwise if $p_{m+1}$ is a composite number then it is removed out before.

Because of the number of elements in a period is $\prod\limits_{i=0}\limits^{m} (p_{i}-1)>0$, $\forall m\geq 0$, hence the sieve procedure is not degenerate, $i.e.$, the sieve steps can be continued indefinitely and we prove that there exist infinitely many primes.

Generalizing the above sieve to other cases, a sieve system is obtained. They can be called the new sieve methods or Liu sieve methods, because of $N=\{1,2,3,...\}$ is obtained from purifying form of [4]. Liu sieve methods operate in the range of all natural numbers, and do not have any error term, they are quite different to the analytic theory of numbers. We call them a part of constructive theory of numbers.

\section{The double sieve}
Take the pre-sieve set to be all natural numbers
\[
N=\{1,2,3,...\}.
\]

Take the sieve set to be all primes $p_{0}=2,p_{1}=3,p_{2}=5,p_{3}=7,p_{4}=11,...$.
Denote the sieve operator $S_{p}$ be a double sieve such that remove all the multiple of prime $p$ and all the multiple of $p$ and $-2$ simultaneously. Since all the multiple of $2$ and all the multiple
of $2$ and $-2$ is the same, $S_{p_{0}}$ is a little special from
another $S_{p}$ ($p\neq2$).

The $0$ step of sieve\\
\centerline{$N_{0} \equiv S_{p_{0}}N = \{1, 3, 5, \cdots \}$.}\\
The ascending period of post-sieve set $N_{0}$ is $2$,~ the number of elements in a period is $1$.

The 1st step of sieve. Remove the multiple of $p_{1} = 3$ and the multiple of $3$ and $-2$ simultaneously we have\\
\centerline{$N_{1} \equiv S_{p_{1}}N_{0} = \{5, 11, 17, 23, 29, \cdots \}$.}\\
The ascending period of post-sieve set $N_{1}$ is $6 = p_{0}p_{1}$, the number of elements in a period is $p_{1}-2 = 1$.

The relation between this step and prime twins is as follows. Prime
twins $p$, $p+2$ with $p>3$ are coprime with $3$. Since one of $p$,
$p+2$, $p+4$ is a multiple of $3$, hence $p+4$ must be multiple of
3. Therefore, when we sieve $3$, $p+4$ and $p+2 = p+4-2$ are removed out
simultaneously. It follows that, only the former number $p$ on prime
twins is the representative element in $N_{1}$. i.e., the later
number $7, 13, 19, 31, \cdots$ of prime twins $5, 7; 11, 13; 17, 19;
29, 31; \cdots$ are removed out by sieve $3$.
The former numbers $5, 11, 17, 29, \cdots$ are remained belong to the post-sieve set $N_{1}$.

The 2nd step of sieve. Remove the multiple of $p_{2} = 5$ and the multiple of $5$ and $-2$ simultaneously we have\\
\centerline{$N_{2} \equiv S_{p_{2}}N_{1} = \{11, 17, 29, 41, 47, 59, \cdots \}$.}\\
The period of $N_{2}$ is $30 = p_{0}p_{1}p_{2}$, the number of elements in a period is $3 = (p_{1}-2)(p_{2}-2)$.

The 3rd step of sieve. Remove the multiple of $p_{3} = 7$ and the multiple of $7$ and $-2$ simultaneously we have
\begin{align*}
  N_{3} \equiv S_{p_{3}}N_{2} &= \{11, 17, 29, 41, 59, 71, 101, 107, 137, 149, \\
   &~~~ 167, 179, 191, 197, 209, 221, \cdots \}.
\end{align*}

The period of $N_{3}$ is $210 = p_{0}p_{1}p_{2}p_{3}$, the number of elements in a period is $15 = (p_{1}-2)(p_{2}-2)(p_{3}-2)$.

Recurrently for the \textit{m}'s step post-sieve set $N_{m}$. Its
period is $\prod^m_{i=0}p_{i}$,~ the number of elements in a period
is $\prod^m_{i=1}(p_{i}-2)$. Since $\prod^m_{i=1}(p_{i}-2)>0
(\forall m \geq 1)$,~ it means that the sieve procedure does not
degenerate,~ so that the sieve step can be continued indefinitely,
and we obtain infinitely many
\noindent post-sieve sets $N_{0}, N_{1}, N_{2}, N_{3}, \cdots$.
In post-sieve sets we arrange the elements in the ascending order (but omit $1$), that is\\
\centerline{$N_0\backslash \{1\}=\{n_{0,1},n_{0,2},\cdots\},$}
\centerline{$N_{m}=\{n_{m,1}, n_{m,2}, \cdots\}$($m=1, 2, \cdots$).}

We call $\{n_{0,1}, n_{1,1}, n_{2,1}, \cdots\}=\{n_{m,1}\}_{m=0}^{\infty}$ the minimum function. Write down the following table.\\

\begin{center}
\begin{tabular}{|c|c|c|c|c|c|c|c|c|c|c|c|c|c|c|c|c|c|}
\hline
$m$     & 0 & 1 & 2 & 3 & 4 & 5 & 6 & 7 & 8 & 9 & 10 & 11 & 12 & 13 & 14 & 15 & 16\\
\hline
$p_m$   & 2 & 3 & 5 & 7 & 11 & 13 & 17 & 19 & 23 & 29 & 31 & 37 & 41 & 43 & 47 & 53 & 59\\
\hline
$n_{m,1}$ & 3 & 5 & 11 & 11 & 17 & 17 & 29 & 29 & 29 & 41 & 41 & 41 & 59 & 59 & 59 & 59 & 71\\
\hline
\end{tabular}
\end{center}
~\\
The above procedure is called the Liu's Sieve Method for prime
twins. And the beginning several terms of $\{n_{m,1},n_{m,1}+2\}_{m=0}^{\infty}$ (drop its repeat part) are
prime twins.

\section{Properties of double sieve}
\begin{lemma}
$\forall m \geq 0$, we have
$n_{m,1}$ and $n_{m,1}+2$ are co-prime with $p_{0},p_{1}, \cdots, p_{m}$.
\end{lemma}
\textbf{Proof.}
$n_{m,1}$ belong to the  post-sieve set of sieve $p_0, p_1,\cdots, p_m$ progressively, $n_{m,1}+2$ belong to the post-sieve set of sieve $p_0$ and $-2$, sieve $p_1$ and $-2, \cdots$, sieve $p_m$ and $-2$ progressively. Hence $n_{m,1}$ and $n_{m,1}+2$ are co-prime with $p_{0},p_{1}, \cdots, p_{m}$.

\begin{lemma}
The minimum function $\{n_{m,1} \}$ is monotone increasing with respect to $m$, but it is not strictly increasing. \\
\end{lemma}

\textbf{Proof.} During we sieve $p_{m}$, the $n_{m-1,1} \in
N_{m}$ has not been removed out, then $n_{m-1,1} = n_{m,1}$,~ this
is the case of $n_{m,1}$ not increasing. During we sieve $p_{m}$,
the $n_{m-1,1}$ has been removed out, then $n_{m-1,1}<n_{m,1}$, this
is the case of $n_{m,1}$ increasing. Combining these two cases, we
prove $lemma ~3.2$.

\begin{definition}
In case $n_{m-1,1}<n_{m,1}$, we call the ordinate number $p_{m}$ the jump point. That is, when we sieve $p_m$, $n_{m-1,1}$ is removed out.
\end{definition}

\begin{lemma}
There exist infinitely many jump points.
\end{lemma}
\textbf{Proof.} If $\{n_{m,1}\}_{m=0}^{\infty}$ has finite jump points, then $\{n_{m,1}\}_{m=0}^{\infty}$ is bounded above. Since $p_{m+1}$ is the smallest element co-prime with $p_{0},p_{1}, \cdots, p_{m}$, therefore by applying $lemma ~3.1$, we have
\[
n_{m,1}\geq p_{m+1} \rightarrow \infty(when~~m \to \infty)~ ,
\]
It is a contradiction. Hence lemma 3.3 is proved.

\section{Effective range of sieve and proof of infinitely many prime twins exist under suitable conditions}
Take the pre-sieve set to be $N=\{1,2,3,\cdots\}$,~ and take the
sieve set to be all primes $\{p_{i}\}_{i=0,1,2,\cdots}$.~ And we
consider the typical case first,~ that is,~ take the sieve operator
to be the single sieve, that is, sieve the multiple of $p_{i}$ only.

Then
\begin{eqnarray*}
  N_{0} \equiv S_{p_{0}}N = \{1, 3, 5, 7, 9, \cdots\},\\
  N_{0} \cap  (p_{0},~ p_{1}^{2}-2] = N_{0} \cap  (2, 7] = \{3, 5, 7\}.
\end{eqnarray*}
$\{3, 5, 7\}$ is called the effective range of sieve $p_{0}=2$.

Similarly,
\begin{eqnarray*}
  N_{1} \equiv S_{p_{1}}N_{0} = \{1, 5, 7, 11, 13, 17, 19, 23, \cdots\},\\
  N_{1} \cap  (p_{1},~ p_{2}^{2}-2] = N_{1} \cap  (3, 23] = \{5, 7, 11, 13, 17, 19, 23\}.
\end{eqnarray*}
$N_{1} \cap  (3, 23]$ is called the effective range of sieve $p_{0}=2$ and $p_{1}=3$.

In general we have $N_{m} \equiv S_{p_{m}}N_{m-1}$.
All the composite numbers of $N_{m}\leq p^2_{m+1}-2$ are removed out when we sieve $p_{0},p_{1},...,p_{m}$, i.e., the effective range of sieve $p_{0},p_{1},...,p_{m}$ are
$$N_{m} \cap  (p_{m},~ p_{m+1}^{2}-2] = \{p_{m+1}, p_{m+2}, \cdots, p_{\tilde{m}}\},$$
where  $\tilde{p}_{m} \leq p_{m+1}^{2}-2$.

Now, we consider the general case,~ i.e. the sieve operator $S_{p_i}$ is the double sieve, that is,
 remove the multiple of $p_{i}$ and the multiple of $p_{i}$ and $-2$ simultaneously.

Then all primes $p$ are divided into $2$ cases.

i. $p+2$ is a composite number,~ $p$ is removed out when $p+2$ is removed out. In this case, we called $p$ a solitary prime.

ii. $p+2$ is a prime. $p$ is removed out only when we sieve $p$.
In this case,~ $p$ is the former number of prime twins.

Therefore when we sieve $p_{0},p_{1}, \cdots, p_{m}$,~ all the composite numbers and solitary primes of
$N_{m} \cap  (p_{m},~ p_{m+1}^{2}-4]$ are sieved out, and only the former numbers of prime twins retain, i.e.,
\begin{eqnarray*}
&&N_{0} \cap  (p_{0},~ p_{1}^{2}-4]  = \{1, 3, 5, 7, 9, \cdots \} \cap (2, 5] = \{3, 5\},\\
&&N_{1} \cap  (p_{1},~ p_{2}^{2}-4] = \{5, 11, 17, 23, 29, \cdots \} \cap (3, 21] = \{5, 11, 17\},\\
&&N_{2} \cap  (p_{2},~ p_{3}^{2}-4]  = \{11, 17, 29, 41, 47, \cdots \} \cap (5, 45] = \{11, 17, 29, 41\},\\
&&N_{3}\cap(p_{3},~ p_{4}^{2}-4] = \{11,17,29,41,59,71,101,107,137,...\}\cap(7,117] \\
  &&~~~~~~~~~~~~~~~~~~~~~~ = \{11,17,29,41,59,71,101,107\}.\\
&& \cdots
\end{eqnarray*}

Denote by $\emptyset$ the empty set. Assume we have the following restrictions:
\begin{equation}
 {\large N_{m}}\cap(p_{m},~ p_{m+1}^{2}-4]\neq\emptyset ~,~ \forall m\geq 0 .
\end{equation}

\begin{lemma}
Under the assumption$(4.1)$ we have
\begin{equation} \label{eqn:4.2}
  n_{m,1}\leq p_{m+1}^{2}-4.
\end{equation}
\end{lemma}
{\bf Proof.} Take element $\ z\in N_{m}\cap (p_{m},~ p_{m+1}^{2}-4]
$,~ then by definition of $\ n_{m,1} $ we have
$$ \large n_{m,1}\leq z \leq p_{m+1}^{2}-4. $$
Lemma 4.1 is proved.

Under the assumption $(4.1)$, we prove that infinitely many prime twins exist.

When $n_{m-1,1}<n_{m,1}$, by definition 3.3, $p_{m}$ is a jump point. It means, when we sieve $p_{m}$, $n_{m-1,1}$ is removed out, hence we have, either
\begin{equation}
  n_{m-1,1}=qp_{m},
\end{equation}
where $q $ is coprime with $\ p_{0},p_{1},...,p_{m-1}$,~ or
\begin{equation}
  n_{m-1,1}+2=\tilde{q}p_{m},
\end{equation}
where $\tilde{q} $ is co-prime with $p_{0},p_{1},...,p_{m-1}$.

Derive from $(4.3)$, we have either
\begin{equation}
  q=1,~~i.e.~~n_{m-1,1}=p_{m},
\end{equation}
or
\begin{equation}
  q\geq p_{m},~~i.e.~~n_{m-1,1}\geq p_{m}^{2}.
\end{equation}

Derive from $(4.4)$, we have either
\begin{equation}
  \tilde{q}=1,~~i.e.~~n_{m-1,1}=p_{m}-2,
\end{equation}
or
\begin{equation}
  \tilde{q}\geq p_{m},~~i.e.~~n_{m-1,1}\geq p_{m}^{2}-2.
\end{equation}

Applying lemma 4.1, we have
\begin{equation}
  n_{m-1,1}\leq p_{m}^{2}-4.
\end{equation}
(4.6) and (4.8) contradict to (4.9).

Since $\ p_{m}\equiv 1(mod~6)$ is removed out by sieve 3, hence $\ p_{m}\equiv -1(mod~6)$, therefore $p_{m}-2\equiv -3(mod~6) $, that is, $p_{m}-2$ is removed out by
sieve $3$, hence $(4.7)$ cannot be true.

Therefore, $(4.5)$ is the only remained case, i.e., $n_{m-1,1}=p_{m}$ is a prime.

Applying lemma 3.1 we have $p_{m}+2=n_{m-1,1}+2$ is co-prime with $p_{0},p_{1},...,p_{m-1}$. And evidently $p_{m}+2$ is co-prime with $p_{m}$. Combining with $(4.9)$ we have $p_{m}+2\leq p_{m}^{2}-2$,
Hence $p_{m}+2 $ cannot be decomposed and it must be  prime and it is $p_{m}+2=p_{m+1}$.  Hence $n_{m-1,1}=p_{m}$ and $p_{m}+2=p_{m+1}$ are prime twins.

The above result combing with lemma 3.3, we have proved the following theorem.
\begin{theorem}
There exist infinitely many prime twins under the assumption $(4.1)$.
\end{theorem}
 We do not know the assumption $(4.1)$ is always valid or not, because when $m$ is large enough, we do not know how many solitary primes between successful $n_{m-1,1}$ and $n_{m,1}$. Perhaps by carefully study of solitary primes we can prove (4.1). Now we still have a question: can we prove theorem 4.2 without assumption $(4.1)$ ? About this question we shall seek some hints by numerical calculation. See the following section.

\section{Numerical calculation for prime twins}
We add difference to the former numbers of prime twins
\begin{center}
$N_{m}\cap (p_{m},~ p_{m+1}^{2}-4],$
\end{center}
such as

$\ N_{0}\cap (p_{0},~ p_{1}^{2}-4]=\{3,5\}$, Now becomes $\ (3,2),5 $.

$\ N_{1}\cap (p_{1},~ p_{2}^{2}-4]=\{5,11,17\}$, Now becomes $\ (5,6),(11,6),17 $.

$\ N_{2}\cap (p_{2},~ p_{3}^{2}-4]=\{11,17,29,41\}$, Now becomes $\ (11,6),(17,12),(29,12),41 $.

$\ N_{3}\cap (p_{3},~ p_{4}^{2}-4]=\{11,17,29,41,59,71,101,107\}$, Now becomes \\
$\ (11,6),(17,12),(29,12),(41,18),(59,12),(71,30),(101,6),107 $.

Write down the $\ m=0,1,..., $ to $\ p_{m}=1129$ $(p_{m+1}^{2}-4=1151^{2}-4=1324797)$
together and drop its repeat part, we obtain a table of 45 pages.
Write down the first and nearly last pages as follows:

\begin{center}
(3,2), (5,6), (11,6), (17,12), (29,12), (41,18), (59,12), (71,30),\\
(101,6), (107,30), (137,12), (149,30), (179,12), (191,6), (197,30),\\
(227,12),(239,30), (269,12), (281,30), (311,36), (347,72), (419,12), \\
(431,30),(461,60), (521,48), (569,30), (599,18), (617,24), (641,18), \\
(659,150),(809,12), (821,6), (827,30), (857,24), (881,138), (1019,12), (1031,18),\\
(1049,12), (1061,30), (1091,60), (1151,78), (1229,48), (1277,12),\\
(1289,12), (1301,18), (1319,108), (1427,24), (1451,30), (1481,6),\\
(1487,120), (1607,12), (1619,48), (1667,30), (1697,24),\\
(1721,66), (1787,84), (1871,6), (1877,54), (1931,18), (1949,48),\\
(1997,30), (2027,54), (2081,6), (2087,24), (2111,18), (2129,12),\\
(2141,96), (2237,30), (2267,42), (2309,30), (2339,42), (2381,168),\\
(2459,42), (2591,66), (2657,30), (2687,24), (2711,18), (2729,60),\\
(2789,12), (2801,168), (2969,30), (2999,120), (3119,48), (3167,84),\\
(3251,6), (3257,42), (3299,30), (3329,30), (3359,12), (3371,18),\\
(3389,72), (3461,6), (3467,60), (3527,12), (3539,18), (3557,24),\\
(3581,90), (3671,96), (3767,54), (3821,30), (3851,66), (3917,12),\\
(3929,72), (4001,18), (4019,30), (4049,42), (4091,36), (4127,30),\\
(4157,60), (4217,12), (4229,12), (4241,18), (4259,12), (4271,66),\\
(4337,84), (4421,60), (4481,36), (4517,30), (4547,90), (4637,12),\\
(4648,72), (4721,66), (4787,12), (4799,132), (4931,36), (4967,42),\\
(5009,12), (5021,78), (5099,132), (5231,48), (5279,138), (5417,24),\\
(5441,36), (5477,24), (5501,18), (5519,120), (5639,12), (5651,6),\\
(5657,84), (5741,108), (5849,18), (5867,12), (5879,210), (6089,42),\\
(6131,66), (6197,72), (6269,30), (6299,60), (6359,90), (6449,102),\\
(6551,18), (6569,90), (6659,30), (6689,12), (6701,60), (6761,18),\\
(6779,12), (6791,36), (6827,42), (6869,78), (6947,12), (6959,168),\\
(7127,84), (7211,96), (7307,24), (7331,18), (7349,108), (7457,30),\\
(7487,60), (7547,12), (7559,30), (7589,168), (7757,120), (7877,72),\\
$\cdots$, \\
(1282907,264), (1283171,366), (1283537,180), (1283717,162),\\
(1283879,60), (1283939,42), (1283981,228), (1284209,342),\\
(1284551,186), (1284737,54), (1284791,258), (1285049,462),\\
(1285511,6), (1285517,30), (1285547,264), (1285811,336),\\
(1286147,42), (1286189,78), (1286267,222), (1286489,330),\\
(1286819,18), (1286837,102), (1286939,42), (1286981,78),\\
(1287059,72), (1287131,66), (1287179,174), (1287371,96),\\
(1287467,84), (1287551,198), (1287749,420), (1288169,78),\\
(1288247,114), (1288361,60), (1288421,120), (1288541,156),\\
(1288697,12), (1288709,120), (1288829,42), (1288871,48),\\
(1288919,678), (1289597,24), (1289621,90), (1289711,36),\\
(1289747,54), (1289801,168), (1289969,198), (1290167,90),\\
(1290257,174), (1290431,36), (1290467,162), (1290629,378),\\
(1291007,12), (1291019,198), (1291217,264), (1291481,336),\\
(1291817,90), (1291907,234), (1292141,450), (1292591,66),\\
(1292657,342), (1292999,318), (1293317,102), (1293419,72),\\
(1293491,348), (1293839,108), (1293947,30), (1293977,42),\\
(1294019,18), (1294037,84), (1294121,78), (1294197,102),\\
(1294301,66), (1294367,282), (1294649,72), (1294721,36),\\
(1294757,312), (1295069,150), (1295219,78), (1295297,24),\\
(1295321,66), (1295387,162), (1295549,12), (1295561,306),\\
(1295867,474), (1296341,180), (1296521,480), (1297001,60),\\
(1297061,108), (1297169,102), (1297271,96), (1297367,30),\\
(1297397,234), (1297631,18), (1297649,462), (1298111,6),\\
(1298117,372), (1298489,162), (1298651,258), (1298909,150),\\
(1299059,150), (1299209,132), (1299341,36), (1299377,60),\\
(1299437,12), 1299449.
\end{center}

From the above table we see that, the successive prime twins
with difference 6 and 12 are both appear frequently.

We study the pair of prime twins with difference 6. The initial pair
is $5,7; 11,13$. The first pair is $11,13; 17, 19$. It is the unique pair in
$(10, 100)$. Then in (100,1000) there have three pairs, they are
$101,103;$ $107,109;$ $191,193;$ $197,199;$ $821,823;$ $827,829$.

Calculating by computer with $p_{m} \leq 9973,~ p^2_{m+1}-4=100140045.$
We obtain the following table for number of pairs in
$(10^{m},10^{m+1})(m=3,4,5,6,7)$.

\begin{center}
\begin{tabular}{c|c}
  ~& \emph{Numbers of pairs for prime twins with difference 6}\\
  \hline
  $(10^3, 10^4)$ & 7\\
  \hline
  $(10^4, 10^5)$ & 26\\
  \hline
  $(10^5, 10^6)$ & 128\\
  \hline
  $(10^6, 10^7)$ & 733\\
  \hline
  $(10^7, 10^8)$ & 3869\\
\end{tabular}
\end{center}

Hence probably we have, number of pairs $\rightarrow\infty$ when $m\rightarrow\infty.$

Take subset to be pairs of prime twins with difference 6, that is $z,z+2$, $w,w+2$ with $w=z+6$ put together, become $\{z,z+2,z+6,z+8\}$, it is a prime quaternary. We show that prime quaternary is a solution of special 4 sieve problem. See the following section.

\section{Special 4 sieve and its properties}
Take the pre-sieve set to be all natural numbers $N=\{1,2,3,4,5,6,...\}$.

Take the sieve set to be all the prime numbers $p_{0}=2,~ p_{1}=3,~ p_{2}=5,~ p_{3}=7,\cdots$. Denote the sieve
operator $S_{p_i}$ be: remove the multiple of $p_i$; remove the multiple of $p_i$ and $-2$; remove the multiple of $p_i$ and $-6$; remove the multiple of $p_i$ and $-8$, simultaneously.

We call the above sieve problem the special $4$ sieve, and study the minimum values except $1$ of post sieve sets.

Evidently the period of pre-sieve set $N$ is $1$.

\textbf{0'th step of sieve.} Denote the post sieve set $S_{p_{0}}N$ to be $N_{0}$, then
\begin{equation*}
    N_{0}=\{1,3,5,7,9,...\}.
\end{equation*}
The period of $N_{0}$ is $2$. The number of elements in a period is $1$.

 \textbf{1st step of sieve.} Denote $S_{p_{1}}N_{0}$ to be $N_{1}$, then
\begin{equation*}
N_{1}=\{5,11,17,23,29,35,41,...\}.
\end{equation*}
The period of $N_{1}$ is $2\times 3=6$. The number of elements in a period is $1$.

\textbf{2nd step of sieve.} Denote $S_{p_{2}}N_{1}$ to be $N_{2}$, then
\begin{equation*}
N_{2}=\{11,41,71,101,131,161,191,...\}.
\end{equation*}
The period of $N_{2}$ is $2\times 3\times 5=30$. The number of elements in a period is $5-4=1$.

\textbf{3rd step of sieve.}  Denote $S_{p_{3}}N_{2}$ to be $N_{3}$, then
\begin{equation*}
N_{3}=\{11,101,191,221,311,401,431,...,821,...,1481...\},
\end{equation*}
where $821=191+3\times 210$ and $1481=11+7\times 210$. The period of $N_{3}$ is $2\times 3\times 5\times 7=210$. The number of elements in a period is $(5-4)(7-4)=3$.

In general in the $m$'th step of sieve. Denote $S_{p_{m}}N_{m-1}$ by $N_{m}$, then the period of $N_{m}$ is $\prod_{i=0}^{m}p_{i}$.

Denote the elements in a period of $N_{m-1}$  by $a_j$, $1\leq j\leq\prod_{i=2}^{m-1}(p_i-4)$, all $a_j$ satisfying the natural restriction $0<a_j<\prod_{i=0}^{m-1}p_i$, where we restrict them in the first period. Then the elements in a period of $N_m$ be $a_j+k\prod_{i=0}^{m-1}p_{i}$, where $0\leq k \leq p_m-1$, $1\leq j\leq \prod_{i=2}^{m-1}(p_i-4)$, and
\begin{equation*}
a_j+k\prod_{i=0}^{m-1}p_{i}\not\equiv 0,-2,-6,-8(mod~p_m).
\end{equation*}

Since $\prod_{i=0}^{m-1}p_{i}$ is co-prime with $p_{m}$, by the theory of linear congruence, the number of elements in a period $\prod_{i=0}^{m}p_{i}$ is $(p_{m}-4)\prod_{i=2}^{m-1}(p_{i}-4)=\prod_{i=2}^{m}(p_{i}-4)$.

Since $\prod_{i=2}^{m}(p_{i}-4)>0, \forall m\geq 2$. Hence the sieve procedure does not degenerate, so that the sieve step can be continued indefinitely. And we obtain infinitely many post sieve sets $N_{m}(m=0,1,2,...)$.

Write down the elements of $N_{m}$ in ascending order, i.e.
\begin{eqnarray*}
  N_{0}\setminus\{1\}&=& \{n_{0,1},n_{0,2},...\}, \\
  N_{m} &=& \{n_{m,1},n_{m,2},...\} \ (m\geq 1).
\end{eqnarray*}

We have

\begin{lemma}
 $n_{m,1},n_{m,1}+2,n_{m,1}+6,n_{m,1}+8$ are co-prime with $p_{0},p_{1},...,p_{m}$, $\forall m \geq 0$.
\end{lemma}
\textbf{Proof.} Like lemma 3.1.

\begin{lemma}
  $n_{m,1}$ is monotone increasing with respect to $m$, but it is not strictly increasing.
\end{lemma}

\textbf{Proof.} Like lemma 3.2.

\begin{definition}
  In case $n_{m-1,1}<n_{m,1}$, we call the ordinate number $p_m$ the jump point.
\end{definition}

\begin{lemma}
  There exist infinitely many jump points.
\end{lemma}

\textbf{Proof.} Like lemma 3.3.

We need a result on single sieve(i.e. when we sieve $p_{i}$, we remove the multiple of $p_{i}$ only) as follows.
\begin{lemma}
   All the composite numbers in $(p_{m},p_{m+1}^{2})$ are removed out by single sieve of $p_{0},p_{1},...,p_{m}$. Hence all elements in post-sieve set are primes.\footnote{It should be noted that the post-sieve set is not empty.}
\end{lemma}
\textbf{Proof.}
Any composite numbers in $(p_{m},p_{m+1}^{2})$ contains a factor less than $p_{m+1}$. And this factor contains a prime factor of $p_{0}, p_{1},...,p_{m}$. Hence lemma 6.4 is proved.

\begin{definition}
   A prime quaternary $\{z,z+2,z+6,z+8\}$ with its first number $z=n_{m,1}~(m\geq 0)$. We call it  a solution of special 4 sieve.
\end{definition}

In 0'th step of sieve $n_{0,1}=3$, $3+2=5$, $3+6=9$, $3+8=11$. Hence we do not get any prime quaternary solution.

In 1st step of sieve $n_{1,1}=p_{2}=5$, $5+2=7$, $5+6=11$, $5+8=13$. Applying lemma 6.1,  $5,7,11,13$ are co-prime with $p_{0},p_{1}$, Applying lemma 6.4,  $5,7,11,13$ are in the prime number part of $(p_{1},p_{2}^{2})=(3,25)$. Hence $5,7,11,13$ are all primes, i.e. $\{5,7,11,13\}$ is a prime quaternary solution of special 4 sieve. We call it the 0'th prime quaternary solution.

Since $n_{2,1}=n_{3,1}=p_4=11$, $11+2=13$, $11+6=17$, $11+8=19$. And $11,13,17,19$ are co-prime with $p_{0},p_{1}$. Hence $11,13,17,19$ are in the prime number part of $(p_{1},p_{2}^{2})=(3,25)$. Hence $\{11,13,17,19\}$ is a prime quaternary solution. We call it the 1st prime quaternary solution.

In 2nd step sieve. Since we sieve $5=p_{2}$, the 0'th prime quaternary solution is removed out. Only the 1st  prime quaternary solution $\{11,13,17,19\}$ retains.

In 3rd step sieve. Except $\{11,13,17,19\}$ retains. We add new $101, 101+2=103,101+6=107,101+8=109$. They  all lie in $(p_{3},p_{4}^{2})=(7,121)$. Applying lemma 6.4,  all composite numbers in $(7,121)$ are sieved out by remove the multiple of $2,3,5,7$. Hence $101,103,107,109$ lie on the prime number part of $(7,121)$. And $\{101,103,107,109\}$ is the new prime quaternary solution of special 4 sieve. We call it the 2nd prime quaternary solution.

In 4th step. We sieve $11$ and $\{11,13,17,19\}$ is removed out. Only the 2nd prime quaternary solution $\{101,103,107,109\}$ retains.

In 5th step. We sieve $p_{5}=13$. Except $\{101,103,107,109\}$ retains. Since $191$ still lies in $N_{5}$. We add new $191,193,197,199$. By lemma 6.1, $191,193,197,199$ are co-prime with $2,3,5,7,11,13$. Hence they lie in the prime number part of $(13,17^{2})=(13,289)$. Applying lemma 6.4, $\{191,193,197,199\}$ is the new prime quaternary solution of special 4 sieve.

When we sieve $2,3,...$ until $23$. Then all the composite numbers in $(23,29^{2})=(23,841)$ are all removed out. But $821$ still lies in $N_{8}$. Hence $\{821,823,827,829\}$ is the new prime quaternary solution.

When we sieve $2,3,...$ until $37$. Then all the composite numbers in $(37,41^{2})=(37,1681)$ are all removes out. But $1481$ still lies in $N_{11}$. Hence $\{1481,1483,1487,1489\}$ is the new prime quaternary solution.

Now we discuss jump points.

From $n_{0,1}=p_{1}=3$, $n_{1,1}=5>3$, hence $3$ is a jump point. The jump point $3$ is not used because of $3,5,9,11$ is not a prime quaternary solution.

From $n_{1,1}=p_{2}=5$, $n_{2,1}>5$, hence $5$ is a jump point. We call $5$ the 0'th jump point, denote by $l_{0}=2,p_{l_{0}}=5$, corresponding to $\{5,7,11,13\}$ is the 0'th prime quaternary solution.

From $n_{2,1}=n_{3,1}=p_{4}=11$, $n_{4,1}>11$, hence $11$ is a jump point. We call $11$ the 1st jump point, denote by $l_{1}=4,p_{l_{1}}=11$, corresponding to $\{11,13,17,19\}$ is the 1st prime quaternary solution.

From $n_{4,1}=n_{5,1}=...=n_{24,1}=p_{25}=101$, $n_{25,1}>101$, hence $101$ is a jump point. We call $101$ the 2nd jump point, denote by $l_{2}=25,p_{l_{2}}=101$, corresponding to $\{101,103,107,109\}$ is the 2nd prime quaternary solution.

From $n_{25,1}=n_{26,1}=...=n_{41,1}=p_{42}=191$, $n_{42,1}>191$, hence $191$ is a jump point. We call $191$ the 3rd jump point, denote by $l_{3}=42,p_{l_{3}}=191$, corresponding to $\{191,193,197,199\}$ is the 3rd prime quaternary solution.

The difference between 0'th jump point $5$ and 1st jump point $11$ is $6$, this situation is special. In general, any jump point greater than $5$ lies in $N_{2}$. Hence, because of the difference between two successive elements in $N_{2}$ is equal to $30$, implies that the difference between two successive jump points $\geq30$.

\section{Biological model for solutions of special 4 sieve problem. And the proof of infinitely many prime twins exist without any restriction}

We don't study the effect range of solution for special 4 sieve, and study the biological model instead.

Solutions of hyperbolic partial differential equations have their
biological model such as life time, life span. Our double sieve
problem, 4 sieve problem, etc. are very like the hyperbolic P.D.E..
Hence, they have biological model also.

Since hyperbolic P.D.E and its biological model do not familiar for all mathematicians, so that we first observe how to distinguish $\{z_i,z_i+2,z_i+6,z_i+8\}$, where $z_i=p_{l_i}$, to be a prime quaternary. A sufficient condition for $\{z_i,z_i+2,z_i+6,z_i+8\}$, where $z_i=p_{l_i}$, to be a prime quaternary are, $z_i,z_i+2,z_i+6,z_i+8$ have no prime factor $\leq (p_{l_i}+8)^{\frac{1}{2}}$. We take it to be no prime factor $< (p_{l_i}+10)^{\frac{1}{2}}$ and denote the largest prime $<(p_{l_i}+10)^{\frac{1}{2}}$ by $p_{h_i}$. We call $p_{h_i}$ and $p_{l_i}$ the birth time and dead time of prime quaternary $\{z_i,z_i+2,z_i+6,z_i+8\}$, where $z_i=p_{l_i}$. And we call time interval $(p_{h_i},p_{l_i})$ the life time or life span of the biological individual $\{z_i,z_i+2,z_i+6,z_i+8\}$, where $z_i=p_{l_i}$, in the biological colony $\{p_{l_i}, p_{l_i}+2, p_{l_i}+6, p_{l_i}+8\}_{i=0}^{\infty}$.

Numerical examples:\\

\begin{center}
\begin{tabular}{|c|c|c|c|c|c|c|c|}
\hline
$i$     & 0 & 1 & 2 & 3 & 4 & 5 & ... \\
\hline
$p_{l_i}$   & 5 & 11 & 101 & 191 & 821 & 1481 & ... \\
\hline
$p_{h_i}$& 3 & 3 & 7 & 13 & 23 & 37 & ... \\
\hline
\end{tabular}
\end{center}

Consequently, we must take $p_{h_i+1}$ to be the smallest prime $\geq (p_{l_i}+10)^{\frac{1}{2}}$. Hence
\begin{eqnarray*}
p_{h_i}<(p_{l_i}+10)^{\frac{1}{2}}\leq p_{h_i+1},\\
p_{h_i}^2-10< p_{l_i}\leq p_{h_i+1}^2-10.
\end{eqnarray*}

In cases $p_{h_i}\geq 5(i\geq2)$, implies $p_{h_i}(p_{h_i}-1)\geq 20$, hence we have

\begin{equation*}
p_{l_i}>p_{h_i}^2-10>p_{h_i}.
\end{equation*}

Therefore

\begin{equation}
p_{l_i}\in (p_{h_i},p_{h_i+1}^2-10].
\end{equation}

Since $p_{h_0}=3$, $p_{l_0}=5$, $p_{h_1}=3$, $p_{l_1}=11$, hence (7.1) is true also for $p_{h_i}<5$.

By definition 6.2, the solution of special 4 sieve is $\{w_i, w_i+2, w_i+6, w_i+8\}$ where $w_i=n_{l_i-1,1}$. We prove the following theorem.

\begin{theorem}
  When $w_i~(\geq5)$ is a jump point of special 4 sieve, then $\{w_i, w_i+2, w_i+6, w_i+8\}$ is a prime quaternary, where $w_{i}=n_{l_{i}-1,1}$.\footnote{In biological model we take the independent variable to be time. So that the jump point is the jump time, it equals to dead time of biological individual.}
\end{theorem}

\textbf{Proof.}
When $p_{l_0}=w_0=5$,  $\{w_0, w_0+2, w_0+6, w_0+8\}=\{5,7,11,13\}$ is a prime quaternary. When $p_{l_i}\geq 11$, we have $p_{l_i}\in (p_{h_i},p_{h_i+1}^2-10]$. Applying lemma 6.1, $w_i, w_i+2, w_i+6, w_i+8$ are co-prime with $p_{0},p_{1},..., p_{l_i-1}$. Moreover, since $w_i$ removes out by sieve $p_{l_i}$, hence we have $w_i=q p_{l_i}$, or $w_i=q p_{l_i}-2$, or $w_i=q p_{l_i}-6$, $w_i=q p_{l_i}-8$$^{**}$, where $q$ is a natural number and $q$ is co-prime with $p_0,p_1,..., p_{l_i-1}$. Therefore we have:
case 1, $q=1$; case 2, $q\geq p_{l_i}$.\footnote{The multipliers of these 4 cases should be denoted by $q_1$, $q_2$, $q_3$, $q_4$. But since we only operate on one multiplier, i.e. any two multipliers do not meet each other. Hence we denote them by a single symbol $q$, which is suitable.}

In case 1, we have $w_i=p_{l_i}$, or $w_i=p_{l_i}-2$, or $w_i=p_{l_i}-6$, or $w_i=p_{l_i}-8$. Since $p_{l_i}>5$ implies $p_{l_i} \in N_2$, it follows that $p_{l_i}\equiv 11(mod~30)$, implies $p_{l_i}\equiv -1(mod~6)$, $p_{l_i}\equiv 1(mod~5)$.

Hence
\begin{eqnarray*}
p_{l_i}-2\equiv -3(mod~6) \equiv 0(mod~3),\\
p_{l_i}-8\equiv -9(mod~6) \equiv 0(mod~3),\\
p_{l_i}-6\equiv -5(mod~5) \equiv 0(mod~5),
\end{eqnarray*}
are all removed out before. Hence we have uniquely $w_i=p_{l_i}\in(p_{h_i},p_{h_i+1}^2-10]$.

In case 2, we have $w_i\geq p_{l_i}^2-8\geq p_{h_i+1}^2-8$ because of $p_{l_i}>p_{h_i}$ implies $p_{l_i}\geq p_{h_i+1}$. Since $w_i=n_{l_i-1,1}$ is the minimum value of case 1 and case 2, therefore combining case 1 and case 2 we have $w_i=p_{l_i}\in(p_{h_i},p_{h_i+1}^2-10]$.

Applying lemma 6.4, all the composite numbers in $(p_{h_i},p^2_{h_i+1})$ are removed out by single sieve of $p_0,p_1,...,p_{h_i}$. Since $\{w_i,w_i+2,w_i+6,w_i+8\}=\{p_{l_i},p_{l_i}+2,p_{l_i}+6,p_{l_i}+8\} \in (p_{h_i},p^2_{h_i+1}-2]$, and $\{w_i,w_i+2,w_i+6,w_i+8\}$ is co-prime with $p_0,p_1,...,p_{h_i}$ because of $l_i-1\geq h_i$. Hence $\{w_i,w_i+2,w_i+6,w_i+8\}$ is a prime quaternary.

The proof of theorem 7.1 is complete.

The inverse of Theorem 7.1 is also true, i.e, when $\{w_i,w_i+2,w_i+6,w_i+8\}$ is a prime quaternary, then $w_i(\geq 5)$  is a jump point of special 4 sieve. The proof is omitted because we do not use it.

Combining theorem 7.1 with lemma 6.3, we obtain that
\begin{theorem}
There exist infinitely many prime quaternary of the following form: $\{w_i,w_i+2,w_i+6,w_i+8\}_{i=0}^{\infty}$ with $w_1-w_0=6$, $w_i-w_{i-1}\geq 30$ $(\forall i>1)$.
\end{theorem}
From theorem 7.2, it follows that
\begin{corollary}
There exist infinitely many prime twins without any restriction.
\end{corollary}

\textbf{Proof.}
Two prime quaternary $5,7,11,13; 11,13,17,19$ implies three prime twins $5,7;11,13;17,19$. Every prime quaternary in $\{w_i,w_i+2,w_i+6,w_i+8\}_{i=2}^{\infty}$  implies two prime twins ${w_i,w_i+2; w_i+6,w_i+8}$.

Hence, Corollary 7.1 is proved.

\section{Appendix}

On pages 52 to 56 of [5], there is a table of  primes less than 5000. From the list we see that there are two prime quaternary between $1$ and $100$: $5,7,11,13$ and $11,13,17,19$, three prime quaternary between $100$ and $1000$: $101,103,107,109$; $191,193,197,199$, and $821,823,827,829$, five prime quaternary between $1000$ and $5000$:
$1481,1483,1487,1489$;
$1871,1873,1877,1879$;
$2081,2083,2087,2089$;
$3251,3253,3257,3259$
and $3461,3463,3467,3469$.
The larger the prime number table is, the more prime quaternary will be obtained. The practical result of [5] is consistent with our theoretical result in section 7.

\section{Acknowledgment}
It is my pleasure to thank Professor Maodong Ye, Yaozhi Cai and Huojun Ruan in department of mathematics, Zhejiang University, for their excellent advice. Especially, I thank Professor Maodong Ye also for he gave me large amount of numerical data related to prime problem calculated by computer. I also thank post doctor of number theory, Deyi Chen in department of mathematics, Zhejiang University, for he read my preprints carefully and also gave many helps.


\begin{thebibliography}{9}
\bibliographystyle{aomalpha}
\bibitem{1} Jingrun Chen, On the representation of a larger even integer as the sum of prime and the product of at most two primes. Sci, Sinica 16(1973), 156-176.
\bibitem{2} Yitang Zhang, Bounded gaps between primes, Annals of Mathematics 179 (2014) (3) 1121-1174.
\bibitem{3} Halhenstam,H., Richert,H.E., Sieve methods., London Math. Soc. Monogr, (1974).
\bibitem{4} Liu Fengsui, The Liu Fengsui's Prime Formula, Internet, Problems and Puzzles, Problems 37(2000).
\bibitem{5} Hua Loo Keng, Introduction to number theory, Springer (1982).
\end{thebibliography}
\end{document}